\documentclass[12pt,twoside]{article}
\usepackage{amsmath,amssymb,amsfonts,theorem,latexsym}

\newtheorem{defn}{Definition}[section]
\newtheorem{lemma}[defn]{Lemma}

	\newtheorem{ex}[defn]{Example}
\newtheorem{thm}[defn]{Theorem}
\newtheorem{prop}[defn]{Proposition}
\newtheorem{cor}[defn]{Corollary}

\newcommand{\h}{{\mathcal H}}

\newcommand{\lk}{\lambda_k}

\newcommand{\ltn}{{\ell}^2(\mathbb N)}

\newcommand{\mn}{\mathbb N}

\newcommand{\mc}{\mathbb C}

\def\bp{{\noindent\bf Proof. \ }}

\def\ep{\hfill$\square$\par\bigskip}

\def\bqs{\begin{equation}}

\def\eqs{\tag*{$\square$}\end{equation}\par\bigskip}

\def\la{\langle}

\def\ra{\rangle}

\def\ftk{\{f_k\}_{k=1}^\infty}

\def\ctk{\{c_k\}_{k=1}^\infty}

\def\gtk{\{g_k \}_{k=1}^\infty}
\def\phitk{\{\phi_k \}_{k=1}^\infty}

\def\etk{\{e_k\}_{k=1}^\infty}

\def\suk{\sum_{k=1}^\infty}

\def\nl{\left|\left|}

\def\nr{\right|\right|}

\def\span{\overline{\text{span}}}

\def\Span{\text{span}}

\def\bop{\begin{op}\rm}

	\def\eop{\end{op}}

\def\du{{\mathcal D}(U)}

\def\bee{\begin{eqnarray}}

\def\ene{\end{eqnarray}}

\def\bes{\begin{eqnarray*}}

	\def\ens{\end{eqnarray*}}

\def\bei{\begin{itemize}}

	\def\eni{\end{itemize}}

\def\lk{\lambda_k}

\def\lj{\lambda_j}

\def\lkn{\{\lambda_k\}_{k=1}^\infty}
\def\htd{H^2(\mathbb{D})}

\def\tnh{\{T^nh\}_{n=0}^\infty}
\def\cjn{\{c_j\}_{j=1}^\infty}
\def\sujj{\sum_{j=1}^\infty}

\def\Tnf{\{T^n f_1\}_{n=0}^\infty}

\begin{document}

\title{\bf\vspace{-39pt} Operator representations of sequences \\ and dynamical sampling}

\author{\vspace{.1in} Ole Christensen, Marzieh Hasannasab, and Diana T. Stoeva 
}

%\date{}

\maketitle

\markboth{\footnotesize \rm \hfill O. Christensen, M. Hasannasab and D.\, T. Stoeva \hfill}
{\footnotesize \rm \hfill Operator representations of sequences and dynamical sampling \hfill}

\abstract{
 This paper is a contribution to the theory of dynamical sampling.
 Our purpose is twofold. We first consider representations
 of sequences in a Hilbert space in terms of iterated actions of a bounded linear operator. This
 generalizes recent results about operator representations of frames, and is motivated by
 the fact that only very special frames have such a representation. As our second contribution we give a new proof of  a construction of a special class of frames that are proved
 by Aldroubi et al. to be representable via a bounded operator. Our proof is based on a
 single result by Shapiro \& Shields and standard frame theory, and our hope is that it eventually can help to
 provide more general classes of frames with such a representation.
\vspace{5mm} 
\vspace{3mm}
}

\section{Introduction}
Dynamical sampling is a new topic in applied harmonic analysis but has already attracted considerable attention \cite{A1,A2,A3,CMPP,FP,AK, olemmaarzieh,olemmaarzieh-E,olemmaarzieh3}.
One of the key questions is how to construct frames $\ftk$
for a separable Hilbert space $\h$ that can be represented on the form $\{T^n\varphi\}_{n=0}^\infty$ for a linear operator $T$ and some $\varphi \in \h.$ It is
known \cite{olemmaarzieh} that
given a sequence $\ftk$ in $\h$ for which $\mbox{span} \ftk$ is infinite-dimensional,
such a representation
is available if and only if $\ftk$ is linearly independent; it is also known that it is
significantly more complicated to obtain such a representation with a bounded operator $T.$
Various characterizations of boundedness have been reported in \cite{olemmaarzieh-E}; prior to
that a class of  such frames were constructed in $\ell^2(\mn)$, based on the
so-called Carleson condition \cite{A1}.

In this paper we extend certain frame results in dynamical sampling to
general sequences. This generalization
is natural for at least two reasons, one of them being the difficulty of
obtaining a representation of a frame in terms of a bounded operator. The second reason is that it turns out that ``nice operator theoretical properties" of $T$
typically is unrelated to frame properties of the underlying sequence; we demonstrate this
point by several concrete examples. In other words: large classes of standard frames in
harmonic analysis are represented by unbounded operators, and sequences that do not form frames might very well have representations in terms of bounded operators.

The second purpose of the paper is to give a detailed analysis of a construction of a class
of frames that can be represented via certain diagonal operators; the construction first
appeared in \cite{A1}. Our  proof is based on just a single result by Shapiro and Shields
\cite{shapiro} and standard
frame theory. We supplement this with a detailed discussion of the Carleson condition, which indeed
is the main ingredient in the construction of such frames.

The paper is organized as follows. Section \ref{70189a} is devoted to the question of
how to obtain a representation of a general sequence in terms of a bounded operator.
This part of the paper does not even use the Hilbert space structure and holds
in general Banach spaces.
The Carleson condition and the associated frames are discussed in Section
\ref{70918b}. 

The standard definitions in frame theory are well-known to the
sampling and signal processing
community, so we will not repeat them here. The only nonstandard terminology
is that we say that a sequence
$\ftk$ in a Hilbert space $\h$ leads to a {\it frame-like expansion} if there
exists a sequence $\gtk$ in $\h$ such that each $f\in \h$ has an expansion of the type
\bee \label{70918c}  f= \suk \la f,g_k\ra f_k.\ene The classical case of a frame-like
expansion is obtained by letting $\ftk$ and $\gtk$ be a pair of dual frames, but the
more general concept of a frame-like expansion also covers other cases, e.g., a pair
of biorthogonal Schauder bases, or a pair of a frame and analysis (or synthesis) pseudo-dual \cite{shidong, Sdualseq}.

\section{Operator representations of sequences}\label{70189a}

In this section we consider linearly independent sequences $\ftk$ in a Hilbert space $\h$
and analyze the existence of a representation $\Tnf$ for a bounded operator 
$T: \mbox{span} \ftk \to \mbox{span} \ftk.$
This generalizes known results for frames, proved in \cite{olemmaarzieh-E}. The generalization is
motivated by the observation
that in general the existence of a representation of
a sequence $\ftk$ of the form $\Tnf$ for a bounded
operator $T$ is not closely related with frame properties of the given sequence $\ftk.$
Let us illustrate this by some examples.

\begin{ex} \label{61105a} Consider a linearly independent frame $\ftk$
for an infinite dimensional Hilbert space $\h$ and the
associated representation $\ftk= \Tnf$ in terms of a linear operator
$T: \mbox{span}\ftk \to \mbox{span}\ftk$. Assume
that $\inf_{k\in\mn} \|f_k\|>0$.    Consider the sequence
$\phitk\subset \h$
given
by $\phi_k:=2^k f_k, \, k\in \mn,$ which leads to a frame-like expansion and satisfies the lower frame condition  but fails the upper one. For any $k\in \mn,$
\bes \phi_{k+1}= 2^{k+1}f_{k+1} =2^{k+1} Tf_k= 2T(2^k f_k)=2T\phi_k.\ens
This shows that $\phitk$ has the representation $\phitk= \{W^n \phi_1\}_{n=0}^\infty$, where $W=2T.$
In particular, the frame $\ftk$ is represented by a bounded operator if and only if
the non-Bessel sequence $\phitk$ is represented by a bounded operator. Consider, e.g.,
the case where $\ftk$ is a Riesz basis for $\h$. Then $\ftk$ is representable by a bounded
operator \cite{olemmaarzieh}, and hence
the  Schauder basis $\phitk$ is representable as $\{W^n \phi_1\}_{n=0}^\infty$ via a bounded operator $W$.
Furthermore, taking $\ftk$ to be a Riesz basis for $\h$, the
family $\{h_k\}_{k=1}^\infty$ given by $h_k:=2^{-k} f_k$ leads to frame-like expansions and it is a Bessel sequence, but does not satisfy the lower frame condition. However, again,
$\{h_k\}_{k=1}^\infty$ is representable as $\{W^n h_1\}_{n=0}^\infty$ via a bounded operator $W$.
\ep \end{ex}

In the following example we show that we can even have a frame-like expansion for a family
of vectors on the form $\{T^n e_1\}_{n=0}^\infty$, where $T$ is bounded and neither the lower
nor the upper frame condition is satisfied.

\begin{ex} \label{61105a} Using that for all $N\in \mn \setminus \{1\},$
\bes 1+ 2 + \cdots + (N-1)=  \frac{(N-1)N}{2}, \, \, \,     1+2 + \cdots + N= \frac{N(N+1)}{2},\ens
and denoting
$$I_N:= \left\{\frac{(N-1)N}{2}, \frac{(N-1)N}{2}+1, \cdots,
\frac{N(N+1)}{2} -1  \right\},$$
we see that
$\mn$ has a splitting  into disjoints sets, $\mn= \bigcup_{N=2}^\infty I_N$. 
Note that $| I_N|=N.$
Let now $\etk$ denote an orthonormal basis for $\h,$ and define the operator $T$ by
\bes Te_k:= \begin{cases} 2e_{k+1}, & k \in I_N, \, N \, \mbox{odd}; \\
	\frac12 e_{k+1},  & k \in I_N, \, N \, \mbox{even}. \end{cases} \ens
Clearly $T$ extends to a bounded linear operator on $\h.$ Furthermore,
\bes \{T_n e_1\}_{n=0}^\infty = \{ e_1, \frac12\, e_2, \frac14\, e_3, \frac12\, e_4, e_5, 2e_6, e_7, \frac12\, e_8, \frac14\, e_9, \frac18\, e_{10}, \dots\}.\ens  Then $\{T_n e_1\}_{n=0}^\infty$ is a Schauder basis, but neither the lower nor the upper
frame condition is satisfied. In fact, $\{T_n e_1\}_{n=0}^\infty$ has a subsequence that tends to infinity in norm,
and another subsequence that tends to zero.
\ep \end{ex}

The following example further confirms that the questions of a sequence $\ftk$ having ``nice frame properties"
and ``a nice representation" $\ftk=\Tnf$ are unrelated in general.  Indeed, the considered family
$\ftk$ is  not a frame and does not even provide a frame-like expansion, but it has a representation $\ftk= \Tnf$ for a bounded and isometric operator $T.$

\begin{ex} \label{exnonexp} Let $\etk$ be an orthonormal basis for a Hilbert space $\h,$ and consider the sequence
$\ftk$ given by
$f_k:=e_k+e_{k+1}, \, k\in \mn.$ By Example 5.4.6 in \cite{CB}, $\ftk$ is a Bessel sequence but not
a frame, despite the fact that $\span \ftk= \h.$ Since $\ftk$ is linearly independent,
we can consider the operator $T: \mbox{span} \ftk \to \mbox{span} \ftk, Tf_k:=f_{k+1},$
and we clearly have $\ftk =   \Tnf.$ Then, for
any $c_k \in \mc,$ and any $N\in \mn,$
\bes \nl T\sum_{k=1}^N c_kf_k  \nr^2 = \nl \sum_{k=1}^N c_k(e_{k+1}+e_{k+2})  \nr^2=
\nl \sum_{k=1}^N c_k(e_{k}+e_{k+1})   \nr^2= \nl \sum_{k=1}^N c_kf_k  \nr^2. \ens
It follows that $T$ has an extension to an isometric operator $T: \h \to \h.$
\ep \end{ex}

Motivated by the above examples we will now consider the question of representability
by a bounded operator for general sequences in a Hilbert space. The reader who checks
the proofs will notice that we do not even use the Hilbert space structure, i.e., the
results in this section hold in Banach spaces as well. Thus, the reader who goes for
the highest generality can think about $\h$ as a separable Banach space and $\la f , g\ra$ as the notation for the action of $g\in \h^{\prime}$ on $f\in \h$.
Given any sequence $\ftk \subset \h,$ define the synthesis operator

\bes U: \du \to \h, \, \, U\ctk:= \sum_{k=1}^\infty c_k f_k,\ens where the domain $\du$ is the set of
all scalar-valued sequences $\ctk$ for which $\sum_{k=1}^\infty c_k f_k$ is convergent. Note
that in contrast to the usual situation in frame analysis, we do not
restrict our attention to sequences $\ctk$ belonging to $\ltn$. Consider the
right-shift operator ${\mathcal T},$ which acts on an arbitrary scalar sequence $\ctk$ by
${\mathcal T}\ctk:= \{0, c_1, c_2, \dots\}.$
A vector space $V$ of scalar-valued sequences
$\ctk$ is said to be {\it invariant under right-shifts} if ${\mathcal T} (V) \subseteq V.$

The following result generalizes one of the key results in
\cite{olemmaarzieh-E}  to the non-frame case.

\sloppy
\begin{thm} \label{61109a} Consider a sequence $\{f_k\}_{k=1}^\infty$
	in $\h$ which has a representation $\{T^n f_1\}_{n=0}^\infty$ for a linear operator $T: \mbox{span} \{f_k\}_{k=1}^\infty \to \mbox{span} \{f_k\}_{k=1}^\infty$. Then the following statements hold.
	\begin{itemize}

			\item[{\rm (i)}]  If $T$ is bounded, then
		the domain $\du$ and the kernel $N_U$ of the
		synthesis operator are invariant under right-shifts, and
		$\{\frac{\|f_{k+1}\|}{\|f_k\|}\}_{k=1}^\infty\in\ell^\infty$ when $f_k\neq 0 $ for all $k\in\mn$.
		
		\item[{\rm (ii)}]
		T is bounded on $ \mbox{span} \{f_k\}_{k=1}^\infty$ if and only if  there is a positive constant $K$ so that
		\begin{equation} \label{eqf}
		\|U {\mathcal T} \{c_k\}_{k=1}^\infty \| \leq K \|U\{c_k\}_{k=1}^\infty \| \ \mbox{for all finite sequences $\{c_k\}_{k=1}^\infty$}.
		\end{equation}

	\end{itemize}
\end{thm}

\bp
Throughout the proof, when $T$ is bounded, we let $\widetilde{T}$ denote its
unique extension to a bounded linear operator on $\span\{f_k\}_{k=1}^\infty$.

(i) Assume that $T$ is bounded and consider first a sequence $\ctk \in \du.$ In order to show that ${\mathcal T}\ctk \in \du,$
i.e., that $\sum_{k=1}^\infty c_k f_{k+1} $ is convergent,
consider any $M,N\in \mn$ with $N>M;$ then
\bes \nl \sum_{k=1}^N c_k f_{k+1} -\sum_{k=1}^M c_k f_{k+1}  \nr &= &\nl \sum_{k=M+1}^N c_k f_{k+1} \nr
= \nl T\sum_{k=M+1}^N c_k f_{k} \nr \\
& \le & ||T|| \,  \nl \sum_{k=M+1}^N c_k f_{k} \nr \to 0 \ \mbox{as $M,N\to \infty$}.\ens 
Thus $\sum_{k=1}^\infty c_k f_{k+1}$ is convergent, i.e., $\du$ is indeed invariant under
right-shifts.

In order to prove the invariance of $N_U$, assume that $\{c_k\}_{k=1}^\infty \in N_U$. The series $\sum_{k=1}^\infty c_k f_{k+1}$ converges by what is already proved, and furthermore
\bes
\sum_{k=1}^\infty c_k f_{k+1} = \sum_{k=1}^\infty c_k Tf_{k}  =\widetilde{T} \sum_{k=1}^\infty c_k f_{k}= 0;\ens
this shows that ${\mathcal T} \{c_k\}_{k=1}^\infty \in N_U,$ as desired.

Finally, for every $k\in\mn$, $\|f_{k+1}\|\leq \|T\|\cdot \|f_k\|$, and thus $\{\frac{\|f_{k+1}\|}{\|f_k\|}\}_{k=1}^\infty\in\ell^\infty$
when $f_k\neq 0$ for all $k\in \mn.$

(ii) Assume first that $T$ is bounded. 
For every $\{c_k\}_{k=1}^\infty\in \du$, we know by (i) that  ${\mathcal T} \{c_k\}_{k=1}^\infty\in\du;$ furthermore, 
$$ \|U {\mathcal T} \{c_k\}_{k=1}^\infty\|= \|\sum_{k=1}^\infty c_{k} f_{k+1}\|
= \|\sum_{k=1}^\infty c_{k} \widetilde{T} f_k\|
\leq 
\|\widetilde{T}\| \cdot\|U\{c_k\}_{k=1}^\infty\|.$$
Clearly this in particular applies to all the finite sequences.

Conversely, assume that there is a constant $K>0$ so that (\ref{eqf}) holds. Take
an arbitrary $f\in \mbox{span}\ftk,$ i.e.,
$f=\sum_{k=1}^N c_k f_k$ for some $N \in \mn$ and some $c_1, \dots, c_N\in \mc;$
letting $c_k=0$ for $k>N$, we have that
\begin{eqnarray*}
	\|Tf\|&=& \|\sum_{k=1}^N c_k f_{k+1}\|= \|\sum_{k=1}^{\infty} c_k f_{k+1}\|
	= \|U {\mathcal T} \{c_k\}_{k=1}^{\infty} \|
	\\
	&\leq& K\|U \{c_k\}_{k=1}^{\infty} \|
	= K \|f\|,
\end{eqnarray*} as desired.
\ep

In Hilbert spaces, the Riesz bases are the ``better bases" compared to general Schauder
bases, due to the unconditional convergence of the frame decomposition. It is also
known that in a Hilbert space the class of Schauder bases that are norm-bounded below and above but do not form Riesz bases, is quite small;
 thus it does not seem to be worthwhile to consider the general Schauder bases in Hilbert
spaces. However, as already stated the results in the current section also
hold if the underlying space is just a Banach space, so it seems natural to
generalize some of the results that are known for Riesz bases in Hilbert spaces
to Schauder bases in Banach spaces.
In Proposition \ref{newpr} we show that for Schauder bases, the invariance of $\du$
under right-shifts in
Theorem \ref{61109a}(i) actually characterizes the case where the
representing operator is bounded. For the special class of  Schauder bases that are scaled Riesz bases,
Proposition \ref{propmriesz} will show that boundedness of $T$ is equivalent to boundedness of a
certain scalar sequence.

Concerning Theorem \ref{61109a}(ii),
as one can see in the proof, the boundedness of $T$ implies validity of the inequality in (\ref{eqf}) not only for the finite sequences, but also for all the elements of $\du$.

Note that when $\ftk$ is a frame
having the form $\Tnf$ for some linear operator $T,$ it is proved in \cite{olemmaarzieh3}   that
$T$ is bounded if and only if the kernel $N_U$ is invariant under right-shifts. 
Under the weaker assumptions in Proposition \ref{61109a},
boundedness of $T$ still implies the invariance of  $N_U$ under right-shifts, but the converse does not hold, as demonstrated by the following example.

\begin{ex} \label{61109b} Let $\etk$ be a Riesz basis for $\h,$ and define $\ftk\subset \h$ by
$f_k:= k! \, e_k.$  Then $\ftk$ satisfies the lower frame condition, but not the upper frame condition.
Defining $Tf_k=f_{k+1},$ we have
$\ftk= \Tnf;$ however,
\bes Te_k = T( \frac1{k!} f_k) = \frac1{k!}f_{k+1}= (k+1) e_k,\ens
which shows that $T$ is unbounded. On the other hand, $N_U= \{0\},$ so the kernel is
invariant under right-shifts.

Similarly, we can construct a sequence $\ftk$ satisfying the upper frame condition
but not the lower frame condition; for example,
\bes \ftk = \{\frac12\, e_1, \frac13 \, e_2, \frac23\, e_3, \frac14\, e_4, \frac34\, e_5, \dots\}.\ens
Representing this sequence on the form $\ftk= \Tnf$ again leads to an unbounded operator $T,$ while
$N_U= \{0\}$ is invariant under right-shifts.

Notice that in these two examples the domain $\du$ is not invariant under right shifts.
\ep \end{ex}

We will now show that for sequences $\ftk$  leading to a frame-like expansion, an extra
assumption leads to a characterization of the possibility of representing $\ftk$ via a
bounded operator. This generalizes a result in \cite{olemmaarzieh-E}.

\begin{prop}\label{prop2}
	Assume that $\{f_k\}_{k=1}^\infty$ leads to a frame-like expansion  as in
	\eqref{70918c}, and  that $\sum_{k=1}^\infty \langle f, g_k\rangle f_{k+1}$ converges for every $f\in\h$.
	Then $\{f_k\}_{k=1}^\infty$ can be represented as $\{T^nf_1\}_{n=0}^\infty$
	for some bounded linear operator $T:\h\to\h$ if and only if
	\begin{equation}\label{eq2}
	f_{j+1}=\sum_{k=1}^\infty \langle f_j, g_k\rangle f_{k+1}, \ \forall j\in\mn.
	\end{equation}
\end{prop}
\bp
If $\{f_k\}_{k=1}^\infty$ can be represented as $\{T^nf_1\}_{n=0}^\infty$
via some bounded linear operator $T:\h\to\h$,
applying $T$ on the frame-like expansion \eqref{70918c} leads to
$Tf=\sum_{k=1}^\infty \langle f, g_k\rangle f_{k+1}$ for every $f\in\h;$ taking $f=f_j$, $j\in\mn$ now proves (\ref{eq2}).

Conversely, assume that (\ref{eq2}) holds.
Consider the operator $T$ defined by
$$ Tf:=\sum_{k=1}^\infty \langle f, g_k\rangle f_{k+1}, \  f\in\h;$$
note that $T$ is well-defined by assumption and bounded by \cite[Lemma 2.3]{SBuncmult}.
Now, by (\ref{eq2}) it follows that $Tf_j=f_{j+1}$, $j\in\mn$, which proves that
$\ftk= \{T^nf_1\}_{n=0}^\infty,$ as desired.
\ep

Notice that  $\sum_{k=1}^\infty \langle f, g_k\rangle f_{k+1}$ may  converge for all $f\in\h$ even in cases where $\{f_k\}_{k=1}^\infty$ and/or $\{g_k\}_{k=1}^\infty$ are not frames:

\begin{ex}
	Let $\etk$ be a Riesz basis for $\h$ and consider a sequence $\{m_k\}_{k=1}^\infty$
	of nonzero complex numbers such that $|\frac{m_{k+1}}{m_k}|\leq C, \forall k\in\mn$, for some positive constant $C$. Letting
	$\{f_k\}_{k=1}^\infty=\{m_k e_k\}_{k=1}^\infty$ and
	$\{g_k\}_{k=1}^\infty=\{\frac{1}{m_k} e_k\}_{k=1}^\infty$, the series  $\sum_{k=1}^\infty \langle f, g_k\rangle f_{k+1}$ converges for all $f\in\h$ and thus the conclusion of Proposition \ref{prop2} holds.
	Note that unless  $\{m_k\}_{k=1}^\infty$ is bounded below and above,  $\{f_k\}_{k=1}^\infty$ and
	$\{g_k\}_{k=1}^\infty$ are not frames. \ep
\end{ex}

Any Riesz basis $\ftk$ has a representation $\ftk= \Tnf$ for a linear operator $T,$ which is necessarily bounded \cite{olemmaarzieh}. By linear independence, any  Schauder basis $\ftk$ also has an operator representation as $\Tnf$,
but not necessarily via a bounded operator, as one can see in Example \ref{61109b}.
For Schauder bases we will now
give a necessary and sufficient condition for the boundedness of $T$.

\begin{prop} \label{newpr}
	Let $\ftk$  be a Schauder basis for $\h$ and consider the linear
	operator $T: \Span \ftk \to \Span \ftk$ such that $\ftk= \Tnf$.
	Then  $T$ is bounded if and only if $\du$ is invariant under right-shifts.
\end{prop}
\bp
First assume that $\du$ is invariant under right-shifts and let $\{g_k\}_{k=1}^\infty$ denote the unique biorthogonal sequence associated with $\ftk$.
Let $f\in\h$.
Then  $f=\sum_{k=1}^\infty \langle f,g_k\rangle f_k$. Therefore $\{\la f,g_k\ra\}_{k=1}^\infty\in\du$ and hence, by the right-shift invariance of $\du$, the series $\sum_{k=1}^\infty \la f,g_k\ra f_{k+1}$ converges. Furthermore, for every $j\in\mn$,
$f_{j+1}=\sum_{k=1}^\infty \la f_j, g_k\ra f_{k+1}$.
It now follows from Proposition \ref{prop2} that $T$ is bounded.
The converse implication is given in Theorem \ref{61109a}.
\ep

Recall that Schauder bases might differ from Riesz bases in two aspects: they might not
be norm-bounded below and above, and they might lead to frame-like expansions that
are conditionally convergent. For the class of Schauder bases consisting of scaled
Riesz bases, we can characterize the  boundedness of $T$ in a way that is much easier to check compared to the right-shift invariance of $\du$:

\begin{prop} \label{propmriesz} 
Let $\{f_k\}_{k=1}^\infty=\{m_k e_k\}_{k=1}^\infty$, where $\{e_k\}_{k=1}^\infty$ is a Riesz basis for $\h$ and $\{m_k\}_{k=1}^\infty$ is a sequence of  nonzero complex numbers, and consider the operator $T:\mbox{span} \{f_k\}_{k=1}^\infty \to \mbox{span} \{f_k\}_{k=1}^\infty$ such that $\{f_k\}_{k=1}^\infty=\{T^nf_1\}_{n=0}^\infty$. Then $T$ is bounded if and only if  $\{\frac{\|f_{k+1}\|}{\|f_k\|}\}_{k=1}^\infty\in\ell^\infty(\mn)$. 
\end{prop}
\bp If $T$ is bounded, the conclusion follows from Theorem \ref{61109a}.
For the converse part, 
 assume that $\{\frac{\|f_{k+1}\|}{\|f_k\|}\}_{k=1}^\infty\in\ell^\infty$. Then there is a positive constant $C$ such that  $\frac{|m_{k+1}|}{|m_k|}\leq C$ for every $k\in\mn$.
Now let $A$ and $B$ denote  Riesz basis bounds of $\{e_k\}_{k=1}^\infty$.
For every finite sequence $\{c_k\}$ we have
\begin{eqnarray*}
	\|T\sum c_k e_k\|^2 &=&\|\sum \frac{c_k}{m_k} f_{k+1}\|^2 =\|\sum c_k\frac{m_{k+1}}{m_k} e_{k+1}\|^2
	\\
	&\leq& B \sum |c_k\frac{m_{k+1}}{m_k}|^2
	\leq B C^2 \sum |c_k|^2
	\leq \frac{B}{A} C^2 \|\sum c_k e_k\|^2;
\end{eqnarray*}
thus the operator $T$ is bouned on $\mbox{span} \{e_k\}_{k=1}^\infty=\mbox{span} \{f_k\}_{k=1}^\infty$.
\ep

We will now show that if a sequence  consists of a Schauder basis and
a finite and strictly positive number of additional elements, then it  can not be
represented by a bounded operator $T$; this extends a result in \cite{olemmaarzieh}.

\begin{prop} Let $\{f_k\}_{k=1}^\infty$ be a linearly independent sequence in $\h$ containing an  $\omega$-independent  subsequence $\{f_k\}_{k=N+1}^\infty$ for some $N\in\mn$.
	Furthermore, assume that at least one $f_{i_0}$, $i_0\in \{1,\ldots,N\}$, can be written as $\sum_{k=N+1}^\infty c_k f_k$
	for some scalar coefficients $c_k\in \mc.$ Then the operator $T: \Span \ftk
	\to \Span \ftk$ such that $\{f_k\}_{k=1}^\infty=\{T^nf_1\}_{n=0}^\infty$ is unbounded.
\end{prop}
\bp 
Assume that $T$ is bounded and extend $T$ by continuity on the closed linear span of $\{f_k\}_{k=1}^\infty$. We split the argument in two cases:

1) If $i_0=N$, then
$ f_{N+1} = Tf_N=T \sum_{k=N+1}^\infty c_k f_k
=\sum_{k=N+1}^\infty c_k f_{k+1},
$
which contradicts   $\{f_k\}_{k=N+1}^\infty$ being $\omega$-independent.

2) If $i_0<N$, then
\begin{eqnarray*}
	f_{i_0+1} &=& Tf_{i_0}=T \sum_{k=N+1}^\infty c_k f_k
	=\sum_{k=N+1}^\infty c_k f_{k+1}\\
	f_{i_0+2} &=& Tf_{i_0+1}=T \sum_{k=N+1}^\infty c_k f_{k+1}
	=\sum_{k=N+1}^\infty c_k f_{k+2}\\
	&\vdots&\\
	f_{N+1} &=&Tf_N=T \sum_{k=N+1}^\infty c_k f_{k+N-i_0}
	=\sum_{k=N+1}^\infty c_k f_{k+N-i_0+1}
\end{eqnarray*}
which contradicts   $\{f_k\}_{k=N+1}^\infty$ being $\omega$-independent.
\ep

\begin{cor}
	Let $\{f_k\}_{k=1}^\infty$ be a linearly independent sequence in $\h$ containing a Schauder basis $\{f_k\}_{k=N+1}^\infty$ for $\h$  for some $N\in\mn$. Then the operator $T: \Span \ftk
	\to \Span \ftk$ such that $\{f_k\}_{k=1}^\infty=\{T^nf_1\}_{n=0}^\infty$ is unbounded.
	
\end{cor}

\section{Frames and the Carleson condition} \label{70918b}
In this section we consider a class of frames which can
be represented via bounded operators. The construction first appeared in
Lemma 3.17 in \cite{A1}; our
purpose is to provide a different proof, which only relies on a single result by Shapiro
and Shields
\cite{shapiro} and standard frame theory. The key ingredient is the so-called Carleson condition on a sequence $\lkn$ in the open unit disc, which we discuss first.

\subsection{The Carleson condition}
Let $\mathbb{D}$ denote the open unit disc in the complex plane. The Hardy space $H^2(\mathbb{D})$ is defined by
\[ \htd:=\left\{ f:\mathbb{D}\to\mathbb{C} \, \big| \, f(z)=\sum_{n=0}^{\infty}a_nz^n,\{a_n \}_{n=0}^\infty\in\ell^2(\mn_0)\right\}.  \]
The Hardy space $H^2(\mathbb{D})$ is a Hilbert space; given
$f,g\in\htd$, $f=\sum_{n=0}^\infty a_nz^n, \\ g=\sum_{n=0}^\infty b_nz^n$, the inner product is  defined by
$$\langle f,g \rangle=\sum_{n=0}^\infty a_n\overline{b_n}.$$
Note that $\{z^n\}_{n=0}^\infty$ is an orthonormal basis for $\htd;$ denoting the canonical
basis for $\ltn$ by $\{\delta_n\}_{n=1}^\infty,$
the operator $\theta:\htd\to\ltn$ defined by $\theta z^n=\delta_{n+1}$ for $n=0,1,\cdots$ is a unitary operator from $\htd$ onto $\ltn$.
\begin{defn}
A sequence $\lkn\subset\mathbb{D}$ satisfies the Carleson condition if
\bee\label{carleson}\displaystyle\inf_{n\in \mn} \prod_{k\neq n}\frac{|\lk-\lambda_n|}{|1-\overline{\lk}\,\lambda_n|}>0.\ene
\end{defn}
Note that at some places in the literature, another terminology is used and a sequence
$\lkn\subset\mathbb{D}$ is said to be {\it uniformly separated } if \eqref{carleson} holds.

For a given sequence $\Lambda=\lkn \subset \mathbb{D}$, define the sequence-valued operator $\Phi_\Lambda$ by
\bee\label{906a} \Phi_\Lambda f=\{f(\lk)\sqrt{1-|\lk|^2}\}_{k=1}^\infty, \quad f\in\htd.\ene Note that the sequence in \eqref{906a} does not necessarily
belong to $\ell^2(\mn).$
%%%%%%%%%%%%%%%%%%%%%%%%%%%%%%%%%%
In \cite{shapiro}, Shapiro and Shields proved the following result.
\begin{prop}\label{706a}
	A sequence $\lkn\subset\mathbb{D}$ satisfies the Carleson condition  if and only if  $\ltn =\Phi_\Lambda H^2(\mathbb{D});$
	in the affirmative case, $\Phi_\Lambda$ is bounded.
\end{prop}

\begin{cor}\label{corcarleson}
The Carleson condition implies that $\suk (1-|\lk|^2)<\infty$ and thus $\lim\limits_{k\to\infty}|\lk|=1$.
\end{cor}
\bp Assume that $\Lambda=\{\lambda_k\}_{k=1}^\infty$ satisfies the Carleson condition. Since the function $f(z)=1, z\in\mathbb{D}$, is in $\htd$, it follows from Proposition \ref{706a} that $\Phi_\Lambda f\in\ltn$. In other words, $\suk (1-|\lk|^2)<\infty$.  \ep

The following result 
(see, e.g., \cite[Thm. 9.2]{duren70}) gives an easily verifiable criterion for the Carleson condition to hold.
\begin{prop}\label{1406a}
	Let $\lkn\subset\mathbb{D}$ be a sequence of distinct numbers. If
	\bee \label{1606c} \exists\, c\in(0,1) \mbox{ such that } \frac{1-|\lambda_{k+1}|}{1-|\lambda_{k}|}\leq c<1, \quad \forall k\in\mn,\ene
	then $\lkn$ satisfies the Carleson condition.
	If $\lkn$ is positive and increasing,  then the condition \eqref{1606c} is also necessary for $\lkn$ to satisfy the Carleson condition.
\end{prop}

%%%%%%%%%%%%%%%%%%%%%%
\begin{cor} For every $\alpha>1$, the sequence $\lkn=\{1-\alpha^{-k}\}_{k=1}^\infty$ satisfies  the Carleson condition.
\end{cor}
\bp Let $\alpha>1$.
The sequence  $\{1-\alpha^{-k}\}_{k=1}^\infty$ is  positive and increasing. Furthermore, for every $k\in\mn$ we have
$$\frac{1-\lambda_{k+1}}{1-\lambda_{k}}=\frac{\alpha^{-k-1}}{\alpha^{-k}}=\frac{1}{\alpha}<1.$$
Thus, by Proposition \ref{1406a}, $\lkn$ satisfies the Carleson condition.
\ep

The following lemma
collects results about modifications on a sequence that preserve the Carleson condition.
\begin{lemma}\label{1506a}
Given any sequence  $\lkn\subset\mathbb{D},$  the following hold:
\bei\item[(i)] If $\lkn$ satisfies  the Carleson condition, then every subsequence of $\lkn$ satisfies the Carleson condition.
\item[(ii)] If $\lk\neq\lj$ for $k\neq j$ and there is some $n\in\mn$ such that $\{\lk\}_{k\geq n}$ satisfies the Carleson condition, then also $\lkn$ satisfies the Carleson condition.
\item[(iii)]If $\lkn$ satisfies  the Carleson condition and $\lk\in[0,1[$
for all $k\in \mn,$ 
then for every $\ell\in\mn$ the sequence $\{\lk^{1/\ell}\}_{k=1}^\infty$ also satisfies  the Carleson condition.
\eni
\end{lemma}
\bp
(i) is straightforward and (ii) is stated in \cite{CMPP}. To prove (iii),  it is enough to show that for any $x,y\in]0,1[$ we have
\bee\label{2205a}
\frac{x^{1/\ell}-y^{1/\ell}}{1-x^{1/\ell}y^{1/\ell}} \geq \frac{x-y}{1-xy}.\ene
Using the identity
$a^\ell-b^\ell 
=(a-b)\sum_{i=0}^{\ell-1}a^{\ell-i-1}b^i,$
we obtain that
\bee\label{2205b}
\frac{x^{1/\ell}-y^{1/\ell}}{1-x^{1/\ell}y^{1/\ell}}=\frac{x-y}{1-xy} f(x,y)
\ene 
where
\bee\label{906c}f(x,y)=\frac{\sum_{i=0}^{\ell-1}(x^{1/\ell}y^{1/\ell})^i}{\sum_{j=0}^{\ell-1}(x^{1/\ell})^{\ell-j-1}(y^{1/\ell})^j}.\ene
Fixing $x\in]0,1[$, a direct calculation shows that  $\frac{\partial}{\partial y}
\displaystyle f(x,y)<0$ for all $y\in ]0,1[$; thus
$f(x,y)$ is decreasing with respect to $y$. Since $\lim\limits_{y\to 1}f(x,y)=1$, we conclude that $f(x,y)\geq 1$ for all $x,y\in]0,1[$.
Using \eqref{2205b}, this proves that \eqref{2205a} holds.
\ep

\subsection{Frame properties and the Carleson condition}
In this section we provide an alternative proof of Lemma 3.19 in \cite{A1}, which yields
a construction of a class of operators $T: \ltn \to \ltn$ for which $\{T^nh\}_{n=0}^\infty$ is a frame for $\ltn$ for certain sequences $h\in \ltn.$ We will formulate the result in the setting of a
general Hilbert space. The original proof of Lemma 3.19 in \cite{A1} uses  properties of the Gramian
associated with a sequence in the underlying Hilbert space, as well as
interpolating sequences in the
Hardy space. We will base our proof on the more
elementary fact that a Bessel sequence is a frame if and only if the frame
operator is surjective, and phrase the interpolation property directly
in terms of surjectivity of the operator $\Phi_\Lambda$ in Proposition \ref{706a}.

Consider a sequence $\lkn\subset\mathbb{D}$ and assume that $\{\sqrt{1-|\lk|^2}\}_{k=1}^\infty\in\ltn$. Given any separable Hilbert space $\h$, choose an orthonormal basis $\etk$ and consider the bounded linear operator $T: \h \to \h$ for which
$Te_k=\lk e_k$. Let $h:=\suk \sqrt{1-|\lk|^2}e_k$ and consider the iterated system
\bee\label{1506b}\{T^nh\}_{n=0}^\infty=\left\{\suk \lk^n\sqrt{1-|\lk|^2}e_k\right\}_{n=0}^\infty.\ene
We will now state the mentioned result from \cite{A1}. Our proof
is only based on  Proposition \ref{706a} and standard frame theory.

\begin{thm}\label{906b} Let $\lkn\subset\mathbb{D}$ and  assume that $\{\sqrt{1-|\lk|^2}\}_{k=1}^\infty\in\ltn$.
	Then the sequence $\{T^nh\}_{n=0}^\infty$ in \eqref{1506b} is a frame for $\h$ if and only if $\lkn$ satisfies the Carleson condition.
\end{thm}
\bp
Define formally the synthesis operator $V:\ell^2(\mn_0)\to\h$ by $V\{a_n\}_{n=0}^\infty=\sum_{n=0}^\infty a_n T^nh$. By Theorem 5.5.1 in \cite{CB}, the sequence $\tnh$ is a frame for $\h$ if and only if the operator $V$ is well-defined and surjective.

First assume that  $\tnh$ is a frame for $\h.$  Take an arbitrary $\cjn\in\ltn$. The surjectivity of $V$ implies that there exists $\{a_n\}_{n=0}^\infty$ in $\ell^2(\mn_0)$ such that $\sum_{n=0}^\infty a_n T^nh=\sujj c_je_j $. It follows that for each $k\in\mn,$
\bee\label{206a}c_k=\langle \sujj c_j e_j,e_k \rangle 
= \sum_{n=0}^\infty a_n \langle  T^nh , e_k \rangle
=
\sum_{n=0}^\infty a_n\lambda_k^n\sqrt{1-|\lambda_k|^2}.\ene
\sloppy
Defining $f\in H^2(\mathbb{D})$ by $f(z)=\sum_{n=0}^{\infty}a_nz^n,$ the equation \eqref{206a} turns into $f(\lambda_k)\sqrt{1-|\lambda_k|^2}=c_k$. Formulated in terms of the operator $\Phi_\Lambda$ in \eqref{906a}, this means that $\ltn\subseteq\Phi_\Lambda H^2(\mathbb{D}) $.
On the other hand, take an arbitrary $f\in H^2(\mathbb{D})$ and choose $\{a_n\}_{n=0}^\infty\in\ell^2(\mn_0)$ such that $f(z)=\sum_{n=0}^{\infty}a_nz^n$. 
For every $k\in\mn$, we have
\bes \la V\{a_n\}_{n=0}^\infty , e_k \ra&=&\la \sum_{n=0}^\infty a_n\sum_{j=1}^\infty \lambda_j^n\sqrt{1-|\lambda_j|^2} e_j , e_k \ra = \sum_{n=0}^\infty a_n\lambda_k^n\sqrt{1-|\lambda_k|^2}. 
\ens 
Therefore, $\Phi_{\Lambda} f=\{\la V\{a_n\}_{n=0}^\infty , e_k \ra\}_{k=1}^\infty\in\ltn$
and hence, $\Phi_\Lambda H^2(\mathbb{D}) \subseteq \ltn$.
Thus we get $\Phi_\Lambda H^2(\mathbb{D}) = \ltn$, which by Proposition \ref{706a} implies that $\lkn$ satisfies the Carleson condition.

Conversely, assume that the sequence $\{\lambda_j\}_{j=1}^\infty \subset {\mathbb D}$ satisfies
the Carleson condition. We first show that $\sum_{n=0}^\infty a_n T^nh$
is convergent for all $\{a_n\}_{n=0}^\infty \in \ell^2(\mn_0).$
Let $\{a_n\}_{n=0}^\infty \in \ell^2(\mn_0),$ and consider the corresponding $f\in H^2({\mathbb D})$ determined by
$f(z):= \sum_{n=0}^\infty a_nz^n.$ By Proposition \ref{706a} we know that $\Phi_\Lambda H^2({\mathbb D})=
\ell^2(\mn),$ so $\{f(\lambda_k) \sqrt{1-|\lambda_k|^2}\}_{k=1}^\infty\in \ell^2(\mn).$ \sloppy
Now for $N\in\mn$, consider the truncated sequence $\{a_n\}_{n=0}^N$
and the associated function $f_N\in H^2({\mathbb D})$ given by
$f_N(z):= \sum_{n=0}^N a_nz^n.$ Again, $\{f_N(\lambda_k) \sqrt{1-|\lambda_k|^2}\}_{k=1}^\infty\in \ell^2(\mn),$ and since $\Phi_\Lambda: H^2({\mathbb D})\to \ell^2(\mn)$  is bounded  by Proposition \ref{706a},
there is a constant $C>0$ such that
\bes || \Phi_\Lambda f- \Phi_\Lambda f_N||^2   \le  C\, ||f-f_N||^2= C\sum_{n=N+1}^\infty |a_n|^2 \to 0 \, \mbox{as} \, N\to \infty.\ens 
It follows that
\bes \sum_{n=0}^N a_n T^nh & = & \sum_{n=0}^N a_n \sum_{k=1}^\infty \sqrt{1-|\lambda_k|^2} \lambda_k^n e_k
=  \sum_{k=1}^\infty \sqrt{1-|\lambda_k|^2}  \sum_{n=0}^N a_n\lambda_k^n e_k
\\ & = & \sum_{k=1}^\infty \sqrt{1-|\lambda_k|^2} f_N(\lambda_k) e_k  \to
\sum_{k=1}^\infty f(\lambda_k) \sqrt{1-|\lambda_k|^2}  e_k \, \mbox{ as } \, N\to \infty. \ens
This proves that $\sum_{n=0}^\infty a_nT^nh$ is convergent  as claimed, and thus $V$ is well defined from $\ell^2(\mn_0)$ into $\h$.
In order to prove that $\{T^nh\}_{n=0}^\infty$ is frame, it is now enough  to show that the synthesis operator $V:\ell^2(\mn_0)\to\h$ 
is surjective. Let $x\in\h$. By Proposition \ref{706a}, there is an $f\in H^2(\mathbb{D})$ such that $f(\lambda_k)\sqrt{1-|\lambda_k|^2}=\langle x,e_k \rangle$ for all $ k\in\mn.$
Choose $\{a_n\}_{n=0}^\infty\in\ell^2(\mn_0)$ such that $f(z)= \sum_{n=0}^\infty a_nz^n.$ Then for each $k\in\mn$,
\bes \langle V\{a_n\}_{n=0}^\infty , e_k \rangle&=& \langle\sum_{n=0}^\infty a_n T^nh , e_k \rangle=
\sum_{n=0}^{\infty}a_n \langle
\sujj \lj^n\sqrt{1-|\lj|^2}e_j , e_k \rangle\\
&=&\sum_{n=0}^{\infty}a_n
\lambda_k^n\sqrt{1-|\lambda_k|^2}=
f(\lambda_k)\sqrt{1-|\lambda_k|^2}=\langle x,e_k \rangle\ens
Therefore $V\{a_n\}_{n=0}^\infty=x$ and thus $V$ is surjective, as desired.
\ep

Note that whenever a sequence $\lkn $ satisfies the Carleson condition, 
the results in Lemma \ref{1506a} and Theorem \ref{906b} immediately lead to a number of alternative frame constructions.

\begin{ex}\label{17057a} Assume that $\lkn\subset ]0,1[$ satisfies the Carleson condition. By Theorem \ref{906b}, the sequence $\{T^nh\}_{n=0}^\infty$ defined in \eqref{1506b} is a frame for $\h$. For a fixed $\ell\in\mn$, define the operator $T_{\ell}$ by $T_{\ell}e_k=\lambda_k^{1/\ell}e_k$, for all $k\in\mn$, and extend it to a bounded operator on $\h$.  Then
\[ \{T_{\ell}^n h\}_{n=0}^\infty=\bigcup_{r=0}^{\ell-1}
\{T_\ell^{n\ell+r}h\}_{n=0}^\infty=\bigcup_{r=0}^{\ell-1} \{T_\ell^{r}T^{n}h\}_{n=0}^\infty. \]
Since each of the operators $T_\ell^{r}$ are surjective and $\{T^n h\}_{n=0}^\infty$ is a frame, it follows that
the sequence $\{T_{\ell}^n h\}_{n=0}^\infty$ is also a frame. 
\ep \end{ex}

\begin{ex} Assume that $\lkn\subset ]0,1[$ satisfies the Carleson condition.
For an arbitrary fixed $\ell\in\mn$, define the operator $T_{\ell}$ by $T_{\ell}e_k=\lambda_k^{1/\ell}e_k$, for all $k\in\mn$, and extend it to a bounded operator on $\h$. By Lemma \ref{1506a}, Corollary \ref{corcarleson}, and Theorem \ref{906b}, $h_{\ell}:=\suk \sqrt{1-|\lambda_k^{1/\ell}|^2}e_k$ is a well-defined element of $\h$ and $\{T_{\ell}^n h_{\ell}\}_{n=0}^\infty$ is a frame for $\h$.
\ep\end{ex}

\vspace{13pt}
\centerline{ACKNOWLEDGEMENT}
\vspace{13pt}
\noindent The authors would like to thank the reviewers for their
comments and suggestions, which improved the presentation.
The third author is grateful for the hospitality of the Technical University of Denmark during her visits there.
The research on this paper was partially supported by the Austrian Science Fund (FWF) START-project FLAME ('Frames and Linear Operators for Acoustical Modeling and Parameter Estimation'; Y 551-N13).
The first two authors thank the Acoustics Research Institute in Vienna for hospitality and
support during visits in 2016 and 2017.

\vspace{.05in}

\noindent
Ole Christensen, Department of Applied Mathematics and Computer Science,  Technical University of  Denmark,  
 Lyngby 2800, Denmark, Email: ochr@dtu.dk\\

\noindent
Marzieh Hasannasab, Department of Applied Mathematics and Computer Science,  Technical University of  Denmark,  
 Lyngby 2800, Denmark, Email: mhas@dtu.dk\\

\noindent
Diana T. Stoeva, Acoustics Research Institute, Austrian Academy of Sciences, 
 Wohllebengasse 12-14, Vienna 1040, Austria, Email: dstoeva@kfs.oeaw.ac.at 

\end{document}